\documentclass[10pt]{article}
\usepackage[english]{babel}										
\usepackage[utf8]{inputenc}										
\usepackage[T1]{fontenc}	
\usepackage{authblk}
\usepackage{float}
\usepackage{amsmath,amsfonts,amssymb,amsthm,cancel,siunitx,
calculator,calc,mathtools,empheq,latexsym}
\usepackage{subfig,epsfig,tikz,float}	
\usepackage{amscd}
\usepackage{booktabs,multicol,multirow,tabularx,array}          

\usepackage{enumitem}
\usepackage{amsmath}
\usepackage{mathrsfs}
\usepackage{rsfso}
\usepackage{hyperref}
\hypersetup{
    colorlinks,
    citecolor=blue,
    filecolor=black,
    linkcolor=black,
    urlcolor=black
}
\usepackage{upgreek}
\usepackage{tikz}
\usetikzlibrary{positioning, braids, decorations.markings}
\definecolor{myblue}{RGB}{80,80,160}
\definecolor{mygreen}{RGB}{50,140,50}
\definecolor{myorange}{RGB}{255,128,0}

\usepackage[backend=biber,
style=alphabetic,
sorting=ynt]{biblatex}
\usepackage{enumitem}
\theoremstyle{definition}

\title{\normalsize\bf%
\uppercase{Brushstrokes and Tensor Products: Painting with a Monoidal Category}
}
\author[1]{Khyathi Komalan \\ kkomalan@caltech.edu}
\affil[1]{California Institute of Technology, United States}

\begin{document}

\date{}

\maketitle

\vspace{-0.5cm}

\bigskip
\noindent
{\small{\bf ABSTRACT.} This article offers an intuitive introduction to monoidal categories through the lens of painting, presenting abstract mathematical concepts with visual and tactile analogies. Aimed at curious undergraduates and non-specialists, it seeks to demystify category theory by showing how ideas like the tensor product, associators, and braidings can be understood as compositional tools on a canvas.

}

\section*{Introduction}

At first glance, painters and category theorists may seem like they inhabit very different worlds — one grounded in colours and brushes, the other in abstraction and arrows. However, in their own unique way, both of them grapple with structure: how things combine, or how local properties can shape global results. This article explores a parallel between the two — specifically, how layering paints and brushstrokes can be modelled through monoidal categories.

We propose an idealized structure where elements of painting — namely, colour, texture, and brushstrokes, along with their associated canvas regions — form a monoidal category. The tensor product represents interactions such as layering or bilateral painting. For simplicity, we restrict this modelling to \textit{alla prima} (wet-on-wet) painting, which lets us abstract away from palette-to-canvas interactions and drying time.

The goal of this model is not to create a physically faithful simulation of painting, but to uncover a structural metaphor: can the logic of brushstrokes and layering — where different logics combine or interfere, or where one action builds on another — be seen as a kind of algebra? And if so, what kind of algebra? By recasting painting as a compositional process, we seek to understand how ideas familiar from category theory — such as the monoidal structure or braiding — can capture the internal logic of painting.

To pursue this metaphor, we introduce a number of simplifying assumptions. The canvas is assumed to be finite, treated as a fixed, discretized set of regions of any finite shape — each a basic unit upon which strokes act. Colors and textures are treated as data rather than physical substances; we do not attempt to simulate pigment absorption or material accumulation. Likewise, brushstrokes are idealized as deterministic morphisms, each encoding a certain effect on a region, such as applying a hue or modifying texture. We ignore physical features like viscosity or drying time, focusing instead on the compositional logic of layering, and focusing on \textit{alla prima} painting allows us to more realistically idealize that a color $A$ can be formed by mixing the right colour(s) with colour $B$.

Within this framework, painting becomes a choreography of morphisms. Brushstrokes compose sequentially, and regions combine via tensor products to model simultaneous or layered actions. First, we define the basic categorical structure of our model: what the objects and morphisms represent, and how the monoidal product captures painterly interaction. Next, we introduce braiding and illustrate it with examples where stroke order matters. Using this setup, we intuitively explain the Yang–Baxter equation for monoidal categories.

This lens, while not a replacement for the experience of painting, offers insight into a structured process — one that can be traced, combined, and reasoned about using the language of category theory. By modelling brushstrokes and textures as morphisms, and treating the canvas as a compositional stage, we uncover a quiet logic behind artistic layering: how gesture accumulates, how order matters, how local actions ripple into global form.

\section*{A Monoidal Structure for Painting}

To begin modelling painting using category theory, we introduce a category we'll call \textbf{Paint}. Think of this as an abstract space where every possible painted state of a canvas — including brushstrokes and textures — is represented. Here, painting becomes a structured system of transformations.

We now define the ingredients of the category: the objects, the morphisms (i.e., arrows), and the monoidal structure that lets us combine them.

First, we define the \textit{objects} of the category, which are "paint states.'' What exactly does it mean to be a "paint state''?

We assume the canvas to be finite, divided into different finite canvas regions of any shape, labelled as $R_{i}$. Each $R_{i}$ is assigned all possible combinations of (colour, texture) that can exist in that region.

For example, in $R_{1}$, we may have (cadmium red, smooth), (eggshell white, stippled), (bright lilac, impasto), and so on, which continues until we have all possible (colour, texture) combinations in $R_{1}$. A single combination like (colour, texture) assigned to a specific region $R_i$ is what we call a \textit{paint state}. These paint states, as well as combinations of them across regions (expressed through the tensor product, which we’ll return to later), form the objects of our category.
\vspace{-12pt}
\begin{figure}[H]
\begin{center}
\includegraphics[width=6cm]{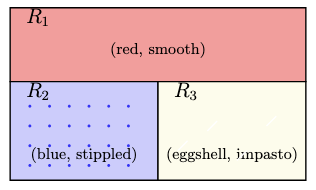}
\caption{Example Paint States Assigned to Canvas Regions}
\end{center}
\end{figure}

Next, we need to define \textit{morphisms} in \textbf{Paint}, which are maps from one object to another. They are defined to be brushstrokes that transform one paint state — i.e., one (color, texture) combination — into another within the same region. Some examples of possible morphisms may be: (a) \verb|add_yellow|: (cadmium red, smooth) in $R_{1}$ $\to$ (orange, smooth) in $R_{1}$; (b) \verb|stipple|: (ultramarine blue, smooth) in $R_{2}$ $\to$ (ultramarine blue, stippled) in $R_{2}$; (c) \verb|blend|: (red, impasto) in $R_{3}$ $\to$ (red-orange, impasto) in $R_{3}$; or (d) \verb|smooth_out|: (sky blue, stippled) in $R_{1}$ $\to$ (sky blue, smooth) in $R_{1}$.
\begin{figure}[H]
    \begin{center}
        \includegraphics[width=8cm]{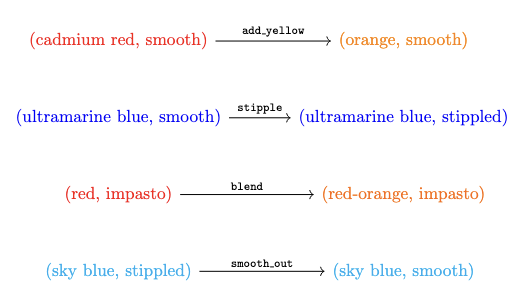}
        \caption{Examples of Morphisms Between Paint States in a Single Region}
    \end{center}
\end{figure}
How do we compose these functions? Let's consider two functions \verb|f|: $(\text{colour}_{1}, \text{texture}_{1})\to (\text{colour}_{2}, \text{texture}_{2})$ and \verb|g|: $(\text{colour}_{2}, \text{texture}_{2})\to (\text{colour}_{3}, \text{texture}_{3}).$ The composition of these functions, \verb|g| $\circ$ \verb|f|, sequentially combines these brush actions, resulting in $(\text{colour}_{1}, \text{texture}_{1})\to (\text{colour}_{3}, \text{texture}_{3}).$ Additionally, all of this occurs in the same region, since all functions act only on the same region.

In a category, there also needs to be an identity morphism, which maps an object to itself. Here, we define it to be \verb|do_nothing|: (colour, texture) $\to$ (colour, texture), which results in no change to the original paint state.

Now that we have a category, we can define a monoidal structure \cite{riehl2017category} on it. 

Perhaps the most essential feature of a \textit{monoidal category} is the \textit{tensor product}, denoted by $\otimes,$ which is \textit{functorial}, but we'll explain what that means later. Simply put, in our case, this tensor product operation takes two painted states on the canvas and produces a new one — a kind of composite that captures how the two inputs interact or layer onto one another. How exactly $\otimes$ behaves depends heavily on the nature of the painted states — most importantly, their location on the canvas. 

To make this concrete, we examine how $\otimes$ behaves in the following cases:

\begin{enumerate}[label=(\alph*)]
    \item \textit{Disjoint regions.}  
    Suppose $(\text{red}, \text{smooth})$ is applied in region $R_1$ and $(\text{blue}, \text{impasto})$ in region $R_2$, with $R_1 \neq R_2$. Then
    \[
    (\text{red}, \text{smooth})_{R_1} \otimes (\text{blue}, \text{impasto})_{R_2}
    \]
    represents a canvas with red smooth paint in $R_1$ and blue impasto in $R_2,$ basically both paint states appear in their respective region in the canvas.

    \item \textit{Same region, different paint states.}  
    Applying $(\text{red}, \text{smooth})$ and $(\text{blue}, \text{transparent})$ both in $R_1$ yields
    \[
    (\text{red}, \text{smooth})_{R_1} \otimes (\text{blue}, \text{transparent})_{R_1},
    \]
    a layered effect — for instance, a purple layer whose appearance depends on how the transparent blue layer interacts optically with the red base.

    \item \textit{Same region, identical states.}  
    Reapplying the same paint state has no further effect:
    \[
    (\text{red}, \text{smooth})_{R_1} \otimes (\text{red}, \text{smooth})_{R_1}
    = (\text{red}, \text{smooth})_{R_1}.
    \]
    That is, $\otimes$ is idempotent when combining identical paint states in the same region.

    \item \textit{Unit object.}  
    The monoidal unit $I$ acts as an identity, $I$ is a blank canvas, and acts on the left and right as follows:
    \[
    I \otimes (\text{red}, \text{smooth})_{R_1}
    = (\text{red}, \text{smooth})_{R_1}
    = (\text{red}, \text{smooth})_{R_1} \otimes I.
    \]

\end{enumerate}

Here's a diagram illustrating the tensor product in each case:

\begin{figure}[H]
    \begin{center}
        \includegraphics[width=10cm]{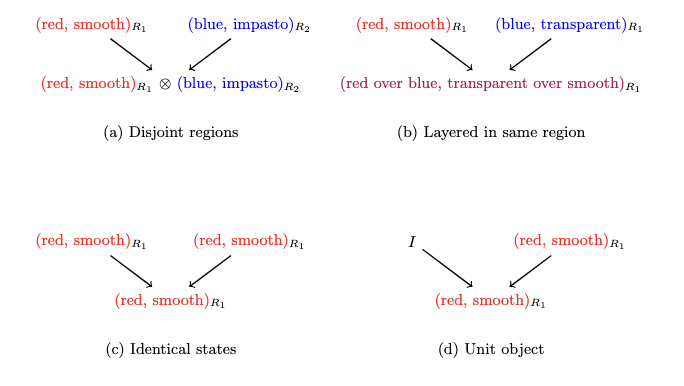}
        \caption{Behavior of the Tensor Product in Paint}
    \end{center}
\end{figure}

Before we move on, we will revisit what exactly it means for our "combining operation'', the tensor product $\otimes,$ to be functorial — some kind of structure preserving map. What exactly does that mean, and why should be care?

Think of it likes this: when we combine two regions on a canvas, we don't just slap them together randomly. If we make changes to a region — like smoothing one out or deepening the colour in the other — we want the final completed painting to reflect those changes faithfully. That's what functoriality ensures. 

Say you have two parts of a painting — one red and smooth, one blue and thickly textured. Now imagine you alter the red part to be darker, and the blue part to be glossier. If you then combine them, functoriality guarantees that the result is exactly what you’d expect: a painting where the red is darker and smooth, and the blue is glossier and textured. It’s what lets painting be a step-by-step process, where changes to each part carry through properly when we combine them.

\begin{figure}[H]
    \begin{center}
        \includegraphics[width=10cm]{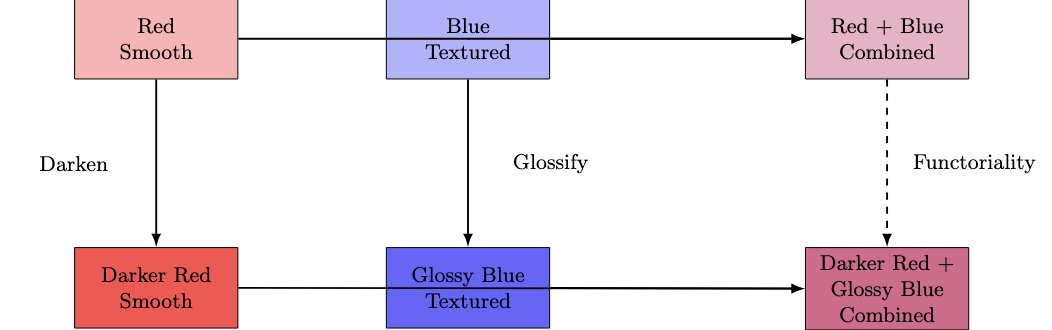}
        \caption{Functoriality in Paint}
    \end{center}
\end{figure}

Just as functoriality ensures that transformed regions behave predictably when combined, something called the \textit{associator} $\alpha$ addresses what happens when we combine three painted regions. Suppose we have three paint states:
\[
(\text{colour}_1, \text{texture}_1)_{R_1},\quad (\text{colour}_2, \text{texture}_2)_{R_2},\quad \text{and} \quad (\text{colour}_3, \text{texture}_3)_{R_3}
\]
There are two natural ways to combine them: we can first blend the paints in \( R_1 \) and \( R_2 \), then apply \( R_3 \) as a layer; or we can first combine \( R_2 \) and \( R_3 \), then layer that mix onto \( R_1 \).

The associator tells us when these combinations are isomorphic to each other. In painting terms, it's like saying: whether we first mix the base and middle layers and then apply the top, or first mix the middle and top layers and then apply that to the base, the final effect is structurally the same. For our specific category, the associator is the identity up to isomorphism, as it doesn't matter which mixing is done first \textit{alla-prima} painting. In other words, we have that:

\[
\alpha: \begin{aligned}
&\big((\text{colour}_1,\, \text{texture}_1)_{R_1} \otimes (\text{colour}_2,\, \text{texture}_2)_{R_2}\big) \otimes (\text{colour}_3,\, \text{texture}_3)_{R_3} \\
&\qquad= (\text{colour}_1,\, \text{texture}_1)_{R_1} \otimes \big((\text{colour}_2,\, \text{texture}_2)_{R_2} \otimes (\text{colour}_3,\, \text{texture}_3)_{R_3}\big)
\end{aligned}
\]

\begin{figure}[H]
    \begin{center}
        \includegraphics[width=10cm, height=8cm]{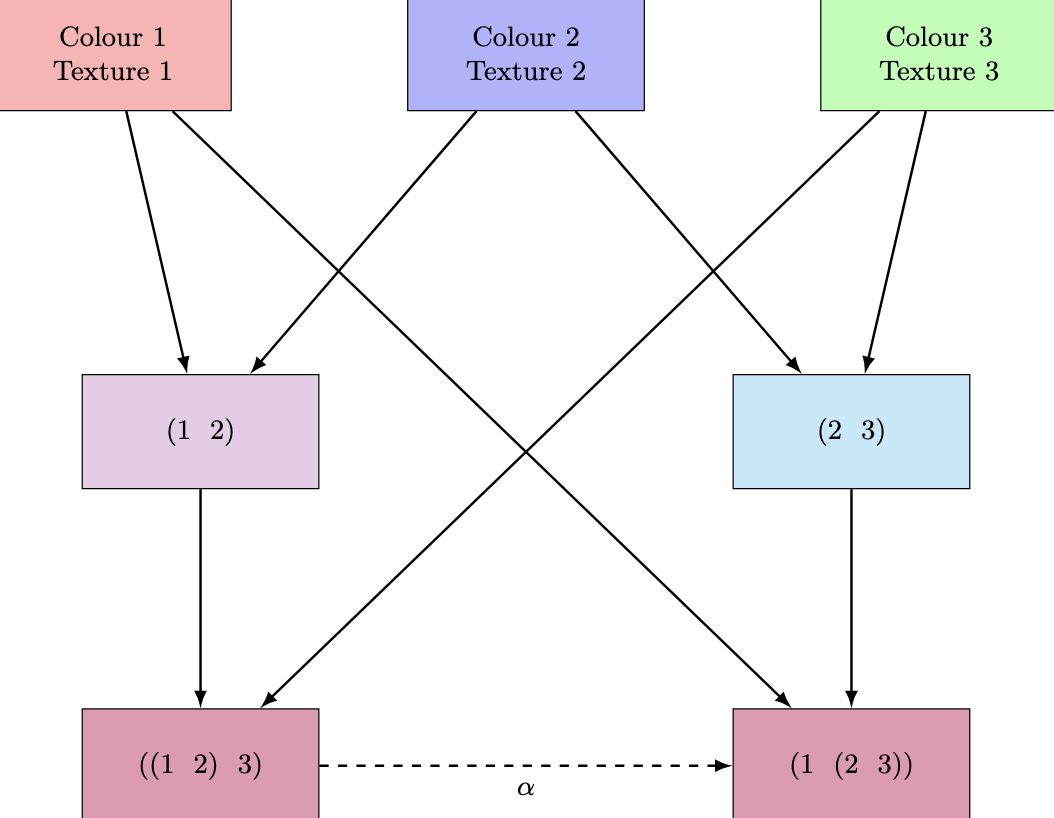}
        \caption{Associator in Paint}
    \end{center}
\end{figure}

Earlier, we introduced the \textit{monoidal unit} — the blank canvas — which acts as both a left and right unit. In general monoidal categories, we have isomorphisms 
\[
\lambda_{\text{object}}:I \otimes \text{object} \cong \text{object} \cong \rho_{\text{object}}:\text{object} \otimes I,
\]
known as the \emph{unitors}. But in our setting, tensoring a painted region with the blank canvas doesn’t change anything: it still results in the same paint state. So, the left and right unitors, $\lambda$ and $\rho$ are identities up to isomorphism.

\begin{figure}[H]
    \begin{center}
        \includegraphics[width=5cm]{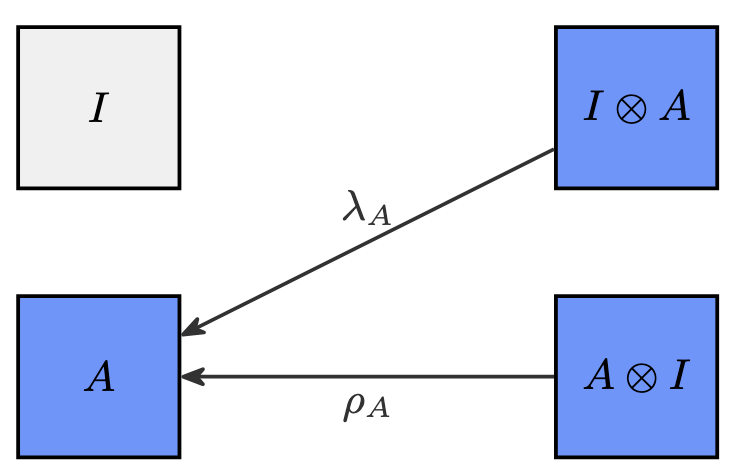}
        \caption{Left and Right Unitors in Paint}
    \end{center}
\end{figure}

Since $\alpha, \lambda$ and $\rho$ are identities, $\textbf{Paint}$ is a \textit{strict} monoidal category.

There are two rules our paint states and tensor products must obey behind the scenes: the \textit{pentagon identity} (which ensures associating four layers of paint in different groupings gives the same result) and the \textit{triangle identity} (which says mixing the blank canvas in doesn’t mess up your workflow). These are backed by giant commuting diagrams shaped like — you guessed it — a pentagon and a triangle.

While discussing the associator, one might think paint can be mixed in any order to produce the same result. This is more true in alla prima painting, where pigments are mixed directly on the canvas. However, the order still matters — especially in how base layers influence saturation. For example, applying red over blue can look quite different from applying blue over red. So, in general,
$$(\text{colour}_1,\, \text{texture}_1)_{R_1} \otimes (\text{colour}_2,\, \text{texture}_2)_{R_2}\neq(\text{colour}_2,\, \text{texture}_2)_{R_2} \otimes (\text{colour}_1,\, \text{texture}_1)_{R_1}$$ 

That said, if the top layer has more pigment (e.g., a heavier brush), it can dominate the mixture — making the result potentially symmetric. This gives us a kind of \textit{commutativity constraint} that depends on physical factors like pigment load.

This context brings us the idea of a \textit{braiding} in a monoidal category \cite{JOYAL199320}. In general, a braiding is a systematic way of switching the order of tensor components while preserving structure — it formalizes a kind of symmetry, even when full commutativity doesn't hold.

The braiding maps:
\[
\beta:
(\text{colour}_1,\, \text{texture}_1)_{R_1} \otimes (\text{colour}_2,\, \text{texture}_2)_{R_2}
\longrightarrow
(\text{colour}_2,\, \text{texture}_2)_{R_2} \otimes (\text{colour}_1,\, \text{texture}_1)_{R_1}
\]
which, in the case of distinct colours and textures in the same region, can be thought of as “reversing” the order in which two paint states are layered — say, brushing $(\text{colour}_1,\, \text{texture}_1)_{R_1}$ over $(\text{colour}_2,\, \text{texture}_2)_{R_2},$ versus brushing $(\text{colour}_2,\, \text{texture}_2)_{R_2}$ over $(\text{colour}_1,\, \text{texture}_1)_{R_1}.$ This swap generally doesn't produce an identical result, but under certain physical conditions — like when the top layer has significantly more pigment or opacity — the outcome can approximate a swap. In such cases, the visual result of $(\text{colour}_1,\, \text{texture}_1)_{R_1}$ over $(\text{colour}_2,\, \text{texture}_2)_{R_2}$ might closely resemble that of $(\text{colour}_2,\, \text{texture}_2)_{R_2}$ over $(\text{colour}_1,\, \text{texture}_1)_{R_1},$ depending on pigment load, brush texture, and opacity. That is, the morphism $\beta$ captures how much influence the second stroke has over the first.

Crucially, this is not full symmetry — paint layering is inherently asymmetric in most physical cases — but it does define a structured transformation. The braiding becomes a controlled, interpretable way of rearranging how strokes are layered, mediated by how dominant one layer is over the other.

For example, if red is brushed lightly over blue, the blue base dominates; but if red is applied with a heavily loaded brush (more pigment, higher opacity), it can override the blue underneath — essentially inverting the layering order's visual effect. In this sense, the braiding morphism is sensitive to pigment dominance, making it a physically meaningful reinterpretation of categorical structure. 

\begin{figure}[H]
    \begin{center}
        \includegraphics[width=9cm]{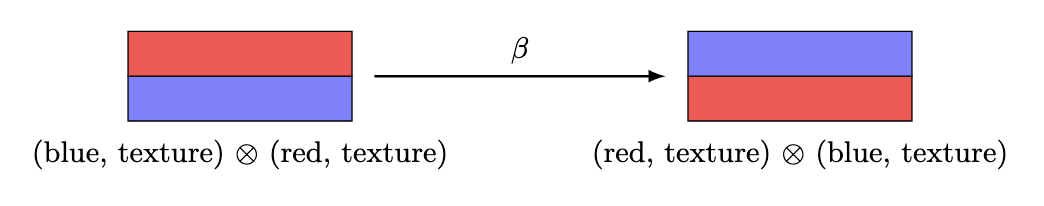}
        \caption{Example of Braiding in Paint}
    \end{center}
\end{figure}

However, in our category, if the two paint states involved in the braiding morphism are located in distinct regions — that is, $R_1 \neq R_2$ — then their interaction is physically disjoint. As a result, swapping their order has no observable effect, and the braiding morphism reduces to the identity, also known as the \textit{trivial braiding}. Similarly, if one of the objects is the unit object $I$ (representing the absence of paint in a region), then braiding has no visual consequence, and again the morphism is trivial. Thus, nontrivial braiding only arises when both paint states occupy the same region and have meaningful physical interaction.

\begin{figure}[H]
    \begin{center}
        \includegraphics[width=9cm]{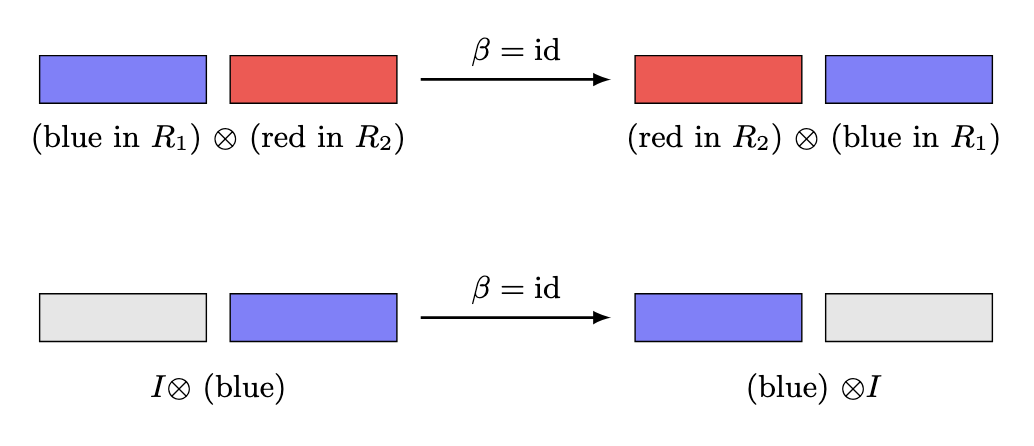}
        \caption{Examples of Trivial Braiding in Paint}
    \end{center}
\end{figure}

There's an elegant consequence that a braiding on a strict monoidal category satisfies — the Yang-Baxter equation \cite{etingof2015tensor}, which shows how three objects interact when braided in sequence.

In painting terms, imagine you have three regions on the canvas — say, a patch of indigo at the top, a smear of vermilion across the middle, and a translucent yellow wash at the bottom. You might drag these colours over one another in different orders: swap the top two, then the bottom two, then the top again; or do the same swaps but in a different order. Intuitively, you’d expect the final composition — the way colours interact and layer — to be the same no matter which sequence you followed, as long as all the same interactions happened. The Yang–Baxter equation is the mathematical guarantee of this: a kind of \textit{consistency of interference}. 

Technically, the Yang-Baxter equation implies that in a strict monoidal category $C$ with braiding $\beta,$ for objects $X, Y, Z \in C,$ the two different ways of "braiding three things'' yield the same result:
$$(\beta_{Y,Z} \otimes \mathrm{id}_X) \circ (\mathrm{id}_Y \otimes \beta_{X,Z}) \circ (\beta_{X,Y} \otimes \mathrm{id}_Z)
=
(\mathrm{id}_Z \otimes \beta_{X,Y}) \circ (\beta_{X,Z} \otimes \mathrm{id}_Y) \circ (\mathrm{id}_X \otimes \beta_{Y,Z})$$

Just as an artist’s brushstrokes may be layered in different sequences but still produce the same overall visual effect, the Yang–Baxter equation ensures that braiding is not just about local swaps, but about \textit{global coherence} in how structure flows. It's a mathematical explanation of how different sequences of mixing paint can ultimately result in the same final mixture.

\section*{What We Learn from Paint}

In modeling painting as a monoidal category, the aim isn’t to impose rigid structure on art or claim that mathematics pervades all things, but to show how engaging with everyday artistic practices can stretch our mathematical imagination. We built a category where objects are painted regions and morphisms are brushstrokes — ways to transform one region into another. We introduced a tensor product to combine paint states, and ensured associativity through an associator that reflected how grouping in composition doesn’t affect the final image. Identity morphisms captured doing nothing — like a brushstroke that leaves no mark — and the unit object represented an empty region. We added a braiding to describe how one stroke can pass over another, and even explored coherence conditions like the Yang-Baxter equation to model consistent layering.

The point wasn’t to reduce painting to math, but to explore how something as intuitive as brushwork can already carry structural assumptions — about order, combination, and interaction — that category theory helps make precise. Thinking this way sharpens our sense of what structure even is, and shows how mathematical ideas can emerge from unlikely places.

\textbf{Funding.} No funding was received to assist with the preparation of this manuscript.

\textbf{Competing Interests.} The author has no relevant financial or non-financial interests to disclose.

\textbf{Data Availability.} No datasets were generated or analysed during the preparation of this manuscript.

\printbibliography

\end{document}